\documentclass[12pt]{article}
\usepackage{amssymb}
\usepackage{amsmath,amsthm}
\usepackage[latin1]{inputenc}
\usepackage{color}
\usepackage{graphicx}

\usepackage{hyperref}
\hypersetup{colorlinks=true, urlcolor=blue}

\usepackage{epsfig}

 \newtheorem{remark}{Remark}

 \newtheorem{theorem}[remark]{Theorem}
 
 \newtheorem{corollary}[remark]{Corollary}

\title{Global defensive $k$-alliances in graphs}

\author{ J. A. Rodr\'{\i}guez-Vel\'{a}zquez$^1$ and J. M. Sigarreta$^2$ \\
\\
$^1${\em Department of Computer Engineering and Mathematics}\\
Rovira i Virgili University of Tarragona \\ Av. Pa\"{\i}sos Catalans
26, 43007 Tarragona, Spain\\
{\small e-mail:\mbox{\tt juanalberto.rodriguez\@@urv.cat}}
\\
$^2${\em Departamento de Matem\'{a}ticas}, Universidad Carlos III de
Madrid,\\ Avda. de la Universidad 30, 28911 Leganés (Madrid)
Spain. \\
{\small e-mail:\mbox{\tt josemaria.sigarreta\@@uc3m.es}} }

\date{}

\begin{document}

\maketitle

\begin{abstract}
Let   $\Gamma=(V,E)$ be a simple graph.  For a nonempty set
$X\subseteq V$, and a vertex $v\in V$, $\delta_{X}(v)$ denotes the
number of
  neighbors  $v$ has in $X$.
A nonempty set $S\subseteq V$ is a \emph{defensive  $k$-alliance} in
$\Gamma=(V,E)$ if $\delta _S(v)\ge \delta_{\bar{S}}(v)+k,$ $\forall
v\in S.$ A defensive $k$-alliance $S$ is called \emph{global} if it
forms a dominating set. The \emph{global defensive $k$-alliance
number} of $\Gamma$, denoted by $\gamma_{k}^{a}(\Gamma)$, is the
minimum cardinality of a defensive $k$-alliance in $\Gamma$. We
study the mathematical properties of $\gamma_{k}^{a}(\Gamma)$.
\end{abstract}



\section{Introduction}

Since (defensive, offensive and dual) alliances  were first
introduced by P. Kristiansen, S. M. Hedetniemi and S. T. Hedetniemi
\cite{alliancesOne}, several authors have studied their mathematical
properties
\cite{cota,chellali,poweful,fava,GlobalalliancesOne,cubic,spectral,planar,kdaf,tesisShafique,line}
as well as the complexity of computing minimum-cardinality of
alliances \cite{np1,np4,np2,np3}. The minimum-cardinality of a
 defensive (respectively, offensive or dual) alliance in a graph $\Gamma$ is called the defensive
 (respectively, offensive or dual)
 alliance number of $\Gamma$. The mathematical properties of defensive
alliances  were first studied in \cite{alliancesOne} where
 several bounds on the defensive
alliance number were given. The particular case of global (strong)
defensive alliances was investigated  in \cite{GlobalalliancesOne}
where several bounds on the global (strong) defensive alliance
number were obtained.  The dual alliances were introduced as
powerful alliances in \cite{cota,poweful}. In \cite{spectral} there
were obtained several tight bounds on the defensive (offensive and
dual) alliance number. In particular, there was investigated the
relationship between the alliance numbers of a graph and its
algebraic connectivity, its spectral radius, and its Laplacian
spectral radius.  Moreover, the study of global defensive (offensive
and dual) alliances in planar graph was initiated in \cite{planar}
and the study of defensive alliances in the line graph of a simple
graph was initiated in cite \cite{line}. The particular case of
global alliances in trees has been investigated in \cite{chellali}.
For many properties of offensive alliances, readers may refer to
\cite{fava,cubic,offensive}.

A generalization of  (defensive and offensive) alliances called
$k$-alliances was presented by K. H. Shafique and R. D. Dutton
\cite{kdaf,kdaf1} where was initiated the study of $k$-alliance free
sets and $k$-alliance cover sets. The aim of this work is to study
mathematical properties of defensive $k$-alliances. We begin by
stating the terminology used. Throughout this article,
$\Gamma=(V,E)$ denotes a simple graph of order $|V|=n$ and size
$|E|=m$. We denote two adjacent vertices $u$ and $v$ by $u\sim v$.
For a nonempty set $X\subseteq V$, and a vertex $v\in V$,
 $N_X(v)$ denotes the set of neighbors  $v$ has in $X$:
$N_X(v):=\{u\in X: u\sim v\},$ and the degree of $v$ in $ X$ will be
denoted by $\delta_{X}(v)=|N_{X}(v)|.$ We denote the degree of a
vertex $v_i\in V$  by $\delta(v_i)$ (or by $d_i$ for short) and the
degree sequence of
 $\Gamma$ by $d_{1}\geq
d_{2}\geq \cdots \geq d_{n}$. The subgraph induced by $S\subset V$
will be denoted by  $\langle S\rangle $ and the complement of the
set $S$ in $V$ will be denoted by $\bar{S}$.

A nonempty set $S\subseteq V$ is a \emph{defensive  $k$-alliance} in
$\Gamma=(V,E)$,  $k\in \{-d_1,\dots,d_1\}$, if for every $ v\in S$,
\begin{equation}\label{cond-A-Defensiva} \delta _S(v)\ge \delta_{\bar{S}}(v)+k.\end{equation}

A vertex $v\in S$ is said to be $k$-\emph{satisfied} by the set $S$
if   (\ref{cond-A-Defensiva}) holds. Notice that
(\ref{cond-A-Defensiva}) is equivalent to
\begin{equation}\label{cond-A-Defensiva1} \delta (v)\ge 2\delta_{\bar{S}}(v)+k.\end{equation}


A defensive $(-1)$-alliance is a \emph{defensive alliance} and a
defensive $0$-alliance is a \emph{strong defensive alliance} as
defined in \cite{alliancesOne}. A defensive $0$-alliance is also
known as a \emph{cohesive set} \cite{porsi}.


The \emph{defensive $k$-alliance number} of $\Gamma$, denoted by
$a_k(\Gamma)$, is defined as the minimum cardinality of a defensive
$k$-alliance in $\Gamma$. Notice that
\begin{equation}a_{k+1}(\Gamma)\ge a_k(\Gamma).\end{equation}

The defensive $(-1)$-alliance number of $\Gamma$ is known as the
\emph{alliance number} of $\Gamma$
 and the defensive  $0$-alliance number is known as  the \emph{strong alliance number},
  \cite{alliancesOne, note, GlobalalliancesOne}.
For instance, in the case of the $3$-cube graph, $\Gamma=Q_{3}$,
every\- set composed by two adjacent vertices is a defensive
alliance of minimum cardinality and every\-  set composed by four
vertices whose induced subgraph is isomorphic to the cycle  $C_4$ is
a strong defensive alliance of minimum cardinality. Thus,
$a_{-1}(Q_3)=2$ and $a_{0}(Q_3)=4$.

For some graphs, there are some values of  $k\in
\{-d_1,\dots,d_1\}$,   such that defensive $k$-alliances do not
exist. For instance, for $k\ge 2$ in the case of the star graph
$S_n$, defensive $k$-alliances do not exist. By
(\ref{cond-A-Defensiva1}) we conclude that, in any graph, there are
defensive $k$-alliances for $k\in \{-d_1,\dots,d_n\}$. For instance,
a defensive $(d_n)$-alliance in $\Gamma=(V,E)$ is $V$. Moreover, if
$v\in V$ is a vertex of minimum degree, $\delta(v)=d_n$, then
$S=\{v\}$ is a defensive $k$-alliance for every $k\le -d_n$.
Therefore, $a_k(\Gamma)=1$, for  $k\le -d_n$.  For the study of the
mathematical properties of $a_k(\Gamma)$, $k\in \{d_n,..., d_1\}$,
we cite \cite{kdefensive}.

A set $S\subset V$ is a  \emph{dominating
set}\label{conjuntodominante} in $\Gamma=(V,E)$ if for every vertex
$u\in \bar{S}$,  $\delta_S(u)>0$ (every vertex in $\bar{S}$ is
adjacent to at least one vertex in S). The \emph{domination number}
of $\Gamma$, denoted by $\gamma(\Gamma)$, is the minimum cardinality
of a dominating set in $\Gamma$.

A defensive $k$-alliance $S$ is called \emph{global} if it forms a
dominating set. The \emph{global defensive $k$-alliance number} of
$\Gamma$, denoted by $\gamma_{k}^{a}(\Gamma)$, is the minimum
cardinality of a defensive $k$-alliance in $\Gamma$. Clearly,
\begin{equation}
\gamma_{k+1}^a(\Gamma)\ge \gamma_k^a(\Gamma)\ge \gamma(\Gamma)\quad
{\rm and } \quad \gamma_k^a(\Gamma)\ge a_k(\Gamma).\end{equation}

The global defensive $(-1)$-alliance number of $\Gamma$ is known as
the \emph{global alliance number} of $\Gamma$
 and the global defensive  $0$-alliance number is known as  the \emph{global strong alliance number}
  \cite{GlobalalliancesOne}. For instance, in the case of the $3$-cube graph, $\Gamma=Q_{3}$,
 every\-  set composed by four
vertices whose induced subgraph is isomorphic to the cycle  $C_4$ is
a global (strong) defensive alliance of minimum cardinality. Thus,
$\gamma^a_{-1}(Q_3)=\gamma^a_{0}(Q_3)=4$.

For some graphs, there are some values of  $k\in
\{-d_1,\dots,d_1\}$,   such that global defensive $k$-alliances do
not exist. For instance, for $k=d_1$ in the case of nonregular
graphs, defensive $k$-alliances do not exist. Therefore, the bounds
showed in this paper on $\gamma^a_k(\Gamma)$, for $ k \le d_1$, are
obtained by supposing that the graph $\Gamma$ contains defensive
$k$-alliances. Notice that for any graph $\Gamma$, every dominating
set is a global defensive $(-d_1)$-alliance. Hence,
$\gamma_{-d_1}^a(\Gamma)=\gamma(\Gamma)$. Moreover, for any
$d_1$-regular graph of order $n$,
$\gamma_{d_1-1}^a(\Gamma)=\gamma_{d_1}^a(\Gamma)=n$.

\section{Global defensive k-alliance number}


\begin{theorem}\label{th1}
Let $S$ be a global defensive $k$-alliance of minimum cardinality in
$\Gamma$. If $W\subset S$ is a dominating set in $\Gamma$, then for
every $r\in \mathbb{Z}$ such that $0 \le r \le
\gamma_k^a(\Gamma)-|W|$,
$$  \gamma_{_{k-2r}}^a(\Gamma)+r\leq \gamma_{_k}^a(\Gamma).$$
\end{theorem}

\begin{proof}
 We can take $X\subset S$ such that $|X|=r.$ Hence, for
every $v\in Y=S-X$,
\begin{align*}
\delta_{Y}(v)&=\delta_S(v)-\delta_X(v)\\
                &\ge \delta_{\bar{S}}(v)+k-\delta_X(v) \\
                &= \delta_{\bar{Y}}(v)+k-2\delta_X(v) \\
                &\ge  \delta_{\bar{Y}}(v)+k-2r.
\end{align*}
Therefore, $Y$ is a defensive $(k-2r)$-alliance in $\Gamma$.
Moreover, as $W\subset Y$, $Y$ is a dominating set and, as a
consequence, $ \gamma_{k-2r}^a(\Gamma)\leq \gamma_k^a(\Gamma)-r.$
\end{proof}

 Notice that if every vertex of $\Gamma$ has even
 degree and $k$ is odd, $k=2l-1$ ,
 then
every defensive $(2l-1)$-alliance in $\Gamma$ is a defensive
$(2l)$-alliance. Hence, in such a case,
$a_{2l-1}(\Gamma)=a_{2l}(\Gamma)$ and
$\gamma^{a}_{2l-1}(\Gamma)=\gamma^{a}_{2l}(\Gamma).$ Analogously, if
every vertex of $\Gamma$ has odd
 degree and $k$ is even, $k=2l$,
 then every defensive $(2l)$-alliance in $\Gamma$ is a defensive
$(2l+1)$-alliance. Hence, in such a case,
$a_{2l}(\Gamma)=a_{2l+1}(\Gamma)$ and
$\gamma^{a}_{2l}(\Gamma)=\gamma^{a}_{2l+1}(\Gamma).$ For instance,
for the complete graph of order $n$ we have
\begin{align*}
n=&\gamma^{a}_{n-1}(K_n)=\gamma^{a}_{n-2}(K_n) \\
 \ge&\gamma^{a}_{n-3}(K_n)=\gamma^{a}_{n-4}(K_n)=n-1\\
&\cdots \\
\ge &\gamma^{a}_{2-n}(K_n)=\gamma^{a}_{3-n}(K_n)=2\\
\ge &\gamma^{a}_{1-n}(K_n)=1.
\end{align*}
Therefore, for every $k\in \{1-n, \dots,  n-1\}$, and for every
$r\in \{0,...,\frac{k+n-1}{2}\}$,
\begin{equation}
\gamma_{_{k-2r}}^a(K_n)+r = \gamma_{_k}^a(K_n).
\end{equation}
Moreover, notice that for every $k\in \{1-n, \dots,  n-1\}$,
 $\gamma_{k}^a(K_n)=\left\lceil
\displaystyle\frac{n+k+1}{2}\right\rceil$.

It was shown in  \cite{GlobalalliancesOne} that
\begin{equation}
 \frac{\sqrt{4n+1}-1}{2}\le \gamma_{-1}^{a}(\Gamma)\le
n-\left\lceil\frac{d_n}{2}\right\rceil \end{equation}
and
\begin{equation}
 \sqrt{n}\le \gamma_0^{a}(\Gamma)\le
n-\left\lfloor\frac{d_n}{2}\right\rfloor.
\end{equation}
Here we generalize the previous results to  defensive $k$-alliances.

\begin{theorem}\label{pilageneral}
For any graph $\Gamma$,
$\displaystyle\frac{\sqrt{4n+k^2}+k}{2} \le \gamma^{a}_k(\Gamma) \le
n-\left\lfloor\displaystyle\frac{d_n-k}{2}\right\rfloor.$
\end{theorem}

\begin{proof}
If $d_n \le k$, then $\gamma^{a}_k(\Gamma) \le n\le
n-\left\lfloor\frac{d_n-k}{2}\right\rfloor.$ Otherwise, consider
 $u\in V$ such that $\delta(u)\ge
\left\lfloor\frac{d_n+d_1}{2}\right\rfloor.$ Let $X\subset V$ be the
set of neighbors  $u$ has in $\Gamma$, $X=\{w\in V: w\sim u\}$. Let
$Y\subset X$ be a vertex set such that
$|Y|=\left\lfloor\frac{d_n-k}{2}\right\rfloor.$ In such a case, the
set $V-Y$ is a global defensive $k$-alliance in $\Gamma$. That is,
$V-Y$ is a dominating set and for every $v\in V-Y$ we have
$\frac{\delta(v)-k}{2}\ge \left\lfloor \frac{d_n-k}{2} \right\rfloor
\ge \delta_{Y}(v)$. Therefore, $\gamma^{a}_k(\Gamma) \le
n-\left\lfloor\frac{d_n-k}{2}\right\rfloor$.

On the other hand, let $S\subseteq V$ be a dominating set in
$\Gamma$. Then,
\begin{equation}\label{cotavecinos3}
n - |S|\le \sum_{v\in S}\delta_{\bar{S}}(v).
\end{equation}
Moreover, if $S$ is a defensive $k$-alliance in $\Gamma$,
 \begin{equation}\label{s(s-1)}
k|S|+\sum_{v\in S}\delta_{\bar{S}}(v)\le \sum_{v\in S}\delta_S(v)\le
|S|(|S|-1).
 \end{equation}
 Hence, solving
 \begin{equation}\label{inecuacion}0\le |S|^2-k|S|-n\end{equation} we deduce the lower bound.
\end{proof}

The upper bound is attained, for instance, for the complete graph
$\Gamma=K_n$ for every $k\in \{1-n, \dots,  n-1\}$.  The lower bound
is attained, for instance, for the 3-cube graph $\Gamma=Q_3$,  in
the following cases: $2\le\gamma_{-3}^a(Q_3)$ and $4\le
\gamma_1(Q_3)=\gamma_0(Q_3)$.

It was shown in  \cite{GlobalalliancesOne} that for any bipartite
graph $\Gamma$ of order $n$ and maximum degree $d_1$,
$$\gamma^{a}_{-1}(\Gamma)\ge \left\lceil
\frac{2n}{d_1+3}\right\rceil \quad {\rm and } \quad
\gamma^{a}_0(\Gamma)\ge \left\lceil \frac{2n}{d_1+2}\right\rceil .$$

Here we generalize the previous bounds to defensive $k$-alliances.
Moreover, we show that the result is not restrictive to the case of
bipartite graphs.

\begin{theorem} \label{nosoloBipartito}
For any graph $\Gamma$,
$\gamma_{k}^{a}(\Gamma)\ge \left\lceil
\displaystyle\frac{n}{\left\lfloor\frac{d_1-k}{2}\right\rfloor
+1}\right\rceil.$
\end{theorem}

\begin{proof}
If $S$ denotes a defensive $k$-alliance in $\Gamma$, then
$$d_1\ge \delta(v)\ge 2\delta_{\bar{S}}(v)+k, \quad \forall v\in
S.$$ Therefore,
\begin{equation} \label{strongGrado}
 \left\lfloor\frac{d_1-k}{2}\right\rfloor \ge \delta_{\bar{S}}(v),
\quad \forall v\in S.
\end{equation}
Hence,
\begin{equation} \label{strongGrado1}
 |S|\left\lfloor\frac{d_1-k}{2}\right\rfloor \ge \sum_{v\in S}\delta_{\bar{S}}(v).
\end{equation} Moreover, if $S$ is a dominating set,  $S$ satisfies
inequality (\ref{cotavecinos3}). The result follows by
(\ref{cotavecinos3}) and (\ref{strongGrado1}).\end{proof}

The above bound is tight. For instance, for the Petersen graph the
bound is attained for every $k$: $3\le\gamma_{-3}^a(\Gamma)$,
$4\le\gamma_{-2}^a(\Gamma)=\gamma_{-1}^a(\Gamma)$, $5\le
\gamma_0(\Gamma)=\gamma_1(\Gamma)$ and $10\le
\gamma_2(\Gamma)=\gamma_3(\Gamma)$. For the 3-cube graph
$\Gamma=Q_3$, the above theorem leads to the following exact values
of $\gamma_{k}^a(Q_3)$: $2\le\gamma_{-3}^a(Q_3)$, $4\le
\gamma_0(Q_3)=\gamma_1(Q_3)$ and $8\le \gamma_2(Q_3)=\gamma_3(Q_3)$.

Hereafter, we denote by ${\cal L}(\Gamma)=(V_l,E_l)$ the line graph
of a simple graph $\Gamma$. The degree of the vertex $e=\{u,v\}\in
V_l$ is $\delta(e)=\delta(u)+\delta(v)-2$. If the degree sequence of
 $\Gamma$ is $d_{1}\geq
d_{2}\geq \cdots \geq d_{n}$, then the maximum degree of ${\cal
L}(\Gamma)$, denoted by  $\Delta_l$, is bounded by
$\Delta_l\leq d_{1}+ d_{2}-2.$

\begin{corollary}
For any graph $\Gamma$ of size $m$ and maximum degrees $d_{1}\ge
d_{2}$,
$$\gamma_{k}^{a}({\cal L}(\Gamma))\geq \left\lceil\frac{m}
{\left\lfloor\frac{d_{1}+d_{2}-2-k}{2}\right\rfloor+1}\right\rceil.$$
 \end{corollary}

The above bound is attained for $k\in \{-3,-2,-1,2,3\}$ in the case
of the complete bipartite graph $\Gamma=K_{1,4}$. Notice that ${\cal
L}(K_{1,4})=K_4$ and $\gamma_{-3}^{a}(K_4)=1$,
$\gamma_{-2}^{a}(K_4)=\gamma_{-1}^{a}(K_4)=2$,
$\gamma_{2}^{a}(K_4)=\gamma_{3}^{a}(K_4)=4$.

In the case of cubic graphs\footnote{A cubic graph is a $3$-regular
graph.} $\gamma(\Gamma)=\gamma_{-3}^a(\Gamma)\le
\gamma_{-2}^a(\Gamma)=\gamma_{-1}^a(\Gamma)\le
\gamma_{0}^a(\Gamma)=\gamma_{1}^a(\Gamma)\le
\gamma_{2}^a(\Gamma)=\gamma_{3}^a(\Gamma)=n$. So, in this case we
only study, $\gamma_{-1}^a(\Gamma)$ and $\gamma_{0}^a(\Gamma)$.

\begin{theorem} \label{lineas}
For any cubic
 graph $\Gamma$,
$\gamma^{a}_{-1}(\Gamma)\leq 2\gamma(\Gamma).$
\end{theorem}

\begin{proof}
Let $S$ be a dominating set of minimum cardinality in $\Gamma$. Let
$X\subseteq S$ be the set composed  by all $v_i\in S$ such that
$\delta_S(v_i)=0$. For each $v_i\in X$ we take a vertex $u_i\in
\bar{S}$ such that $u_i\sim v_i$. Let $Y\subseteq \bar{S}$ defined
as $Y=\displaystyle\cup_{v_i\in X}\{u_i\}$.  Then we have $|Y|\leq
\gamma (\Gamma)$ and the set $S\cup Y$ is a global defensive
(-1)-alliance in $\Gamma$.
\end{proof}

The above bound is tight. For instance, in the case of the 3-cube
graph we have $\gamma^{a}_{-1}(Q_3)=2\gamma(Q_3)=4.$


A set $S\subset V$ is a total dominating set if every vertex in $V$
has a neighbor in $S$. The  total domination number
$\gamma_t(\Gamma)$ is the minimum cardinality of a total dominating
set in $\Gamma$. Notice that if $\Gamma$ is a cubic graph, then
\begin{equation}\gamma^{a}_{-1}(\Gamma)=\gamma_{t}(\Gamma).\end{equation}
It was shown in  \cite{total} that if  $\Gamma$ is a connected graph
of order $n\geq 3$, then \begin{equation}\gamma_{t}(\Gamma)\leq
\frac{2n}{3}.\end{equation} Moreover, by Theorem
\ref{nosoloBipartito} we have \begin{equation}\frac{n}{3}\le
\gamma^{a}_{-1}(\Gamma)\quad {\rm and } \quad \frac{n}{2}\le
\gamma^{a}_{0}(\Gamma).\end{equation}


\section{Defensive $k$-alliances in planar graphs}


It is well-known that the size of a planar graph $\Gamma$ of  order
$n\ge 3$ is bounded by  $m\leq 3(n-2)$. Moreover, in the case of
triangle-free graphs  $m\leq 2(n-2)$. This inequalities allow us
 to obtain tight bounds for the studied parameters.

\begin{theorem}\label{theoremDef}
Let $\Gamma=(V,E)$ be a graph of order $n$.  If $\Gamma$ has a
global  defensive $k$-alliance  $S$ such that the subgraph $\langle
S\rangle$ is planar.
\begin{itemize}
\item[{\rm (}i{\rm )}]{If $n>2(2-k)$, then $|S|\ge
\left\lceil\frac{n+12}{7-k}\right\rceil.$}

\item[{\rm (}ii{\rm )}]{If $n>2(2-k)$ and  $\langle S \rangle$ is a triangle-free graph, then $|S|\ge
\left\lceil\frac{n+8}{5-k}\right\rceil.$}
\end{itemize}
\end{theorem}

\begin{proof}
\mbox{ }
\begin{itemize}
\item[{\rm (}i{\rm )}]{ If $|S|\le 2$, for every $v\in S$ we have
 $\delta_{\bar{S}}(v)\le 1-k$. Thus, $n\le 2(2-k)$. Therefore,
$n>2(2-k) \Rightarrow |S|>2$.

If $\langle S\rangle $  is planar and $|S|>2$, the size of $\langle
S\rangle $ is bounded by
\begin{equation}\label{proofplanar}
\frac{1}{2}\sum_{v\in S}\delta_{S}(v)\le 3(|S|-2).
\end{equation}
If $S$ is a global defensive $k$-alliance in $\Gamma$,
\begin{equation}\label{cotasuma}
k|S|+(n-|S|) \le k|S|+\sum_{v\in S}\delta_{\bar{S}}(v)\le \sum_{v\in
S}\delta_{S}(v).
\end{equation}
By (\ref{proofplanar}) and (\ref{cotasuma}) the result follows.}

\item[{\rm (}ii{\rm )}]{If $\langle S \rangle$ is a triangle-free graph,
then
\begin{equation}\label{proofplanar-Triang-Free}
\frac{1}{2}\sum_{v\in S}\delta_{S}(v)\le 2(|S|-2).
\end{equation}
The result follows by (\ref{cotasuma}) and
(\ref{proofplanar-Triang-Free}).}\end{itemize}\end{proof}

\begin{corollary}\label{corollaryplanar}
  For any planar graph $\Gamma$ of order $n$.
\begin{itemize}
\item[{\rm (a)}]{If $n>2(2-k)$, then
$\gamma_{k}^{a}(\Gamma)\ge
\left\lceil\frac{n+12}{7-k}\right\rceil.$}

 \item[{\rm (b)}]{If
$n>2(2-k)$ and $\Gamma$ is a triangle-free graph, then
$\gamma_{k}^{a}(\Gamma)\ge \left\lceil\frac{n+8}{5-k}\right\rceil.$}
\end{itemize}
\end{corollary}

\begin{figure}[h]
\begin{center}
\caption{ } \label{planarnue}
\includegraphics[angle=90, width=6cm]{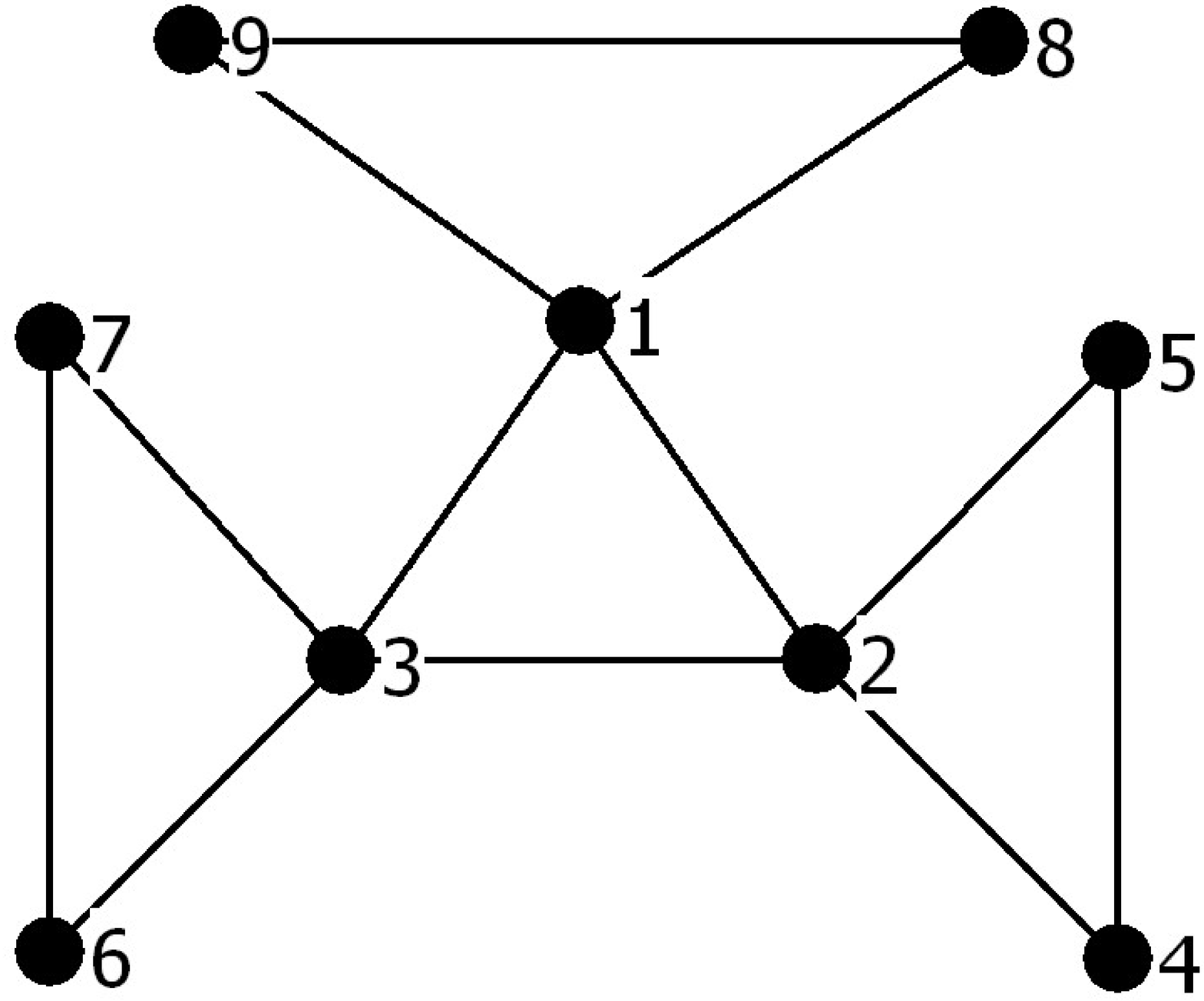}
\end{center}
\end{figure}

 The above bounds are tight. In the case of the  graph of  Figure \ref{planarnue}, the set
$S=\{1,2,3\}$ is a global defensive $k$-alliance for  $k=-2$, $k=-1$
and $k=0$, and Corollary \ref{corollaryplanar}-(a) leads to
$\gamma_{k}^{a}(\Gamma)\geq 3.$ Moreover, if $\Gamma=Q_{3}$, the
$3$-cube graph, Corollary \ref{corollaryplanar}-(b) leads to the
following exact values of $\gamma_{k}^{a}(Q_{3})$: $2\leq
\gamma_{-3}^{a}(Q_{3})$ ,
$4\leq\gamma_{0}^{a}(Q_{3})=\gamma_{1}^{a}(Q_{3})$ and $8 \leq
\gamma_{3}^{a}(Q_{3}).$

\begin{theorem}\label{theorem-faces}
Let $\Gamma$ be a graph of order $n$.  If $\Gamma$ has a global
defensive $k$-alliance  $S$ such that the subgraph $\langle
S\rangle$ is planar connected with $f$ faces. Then,
 $$|S|\ge
\left\lceil\frac{n-2f+4}{3-k}\right\rceil.$$
\end{theorem}

\begin{proof}
By  Euler's formula,  $\displaystyle\sum_{v\in
S}\delta_{S}(v)=2(|S|+f-2)$, and (\ref{cotasuma}) we deduce the
result.\end{proof}

In the case of the  graph of  Figure \ref{planarnue}, the set
$S=\{1,2,3\}$ is a global defensive $k$-alliance for $k=-1$, $k=0$
and $k=2$. Moreover, $\langle S\rangle$ has two faces. In such a
case, Theorem \ref{theorem-faces} leads to $|S|\ge 3$.



\subsection{Defensive $k$-alliances in trees}


In this section we study global defensive $k$-alliances in trees but
we impose a condition on the number of connected components of the
subgraphs induced by the alliances.

\begin{theorem}\label{theoremtree}
Let  $T$ be a tree of order $n$. Let $S$ be a  global defensive
$k$-alliance in $T$ such that the subgraph $\langle S\rangle$
 has $c$ connected components. Then,   $$|S|\ge
\left\lceil\frac{n+2c}{3-k}\right\rceil.$$
\end{theorem}

\begin{proof}
As the subgraph $\langle S\rangle$ is a forest with $c$ connected
components,
\begin{equation}\label{prooftree}
\sum_{v\in S}\delta_S(v)=2(|S|-c).
\end{equation}
The bound of $|S|$ follows from (\ref{cotasuma}) and
(\ref{prooftree}).\end{proof}

\begin{figure}[h]
\begin{center}
 \caption{ } \label{ej3}
 \vspace{0,5cm}
\includegraphics[angle=90, width=6cm]{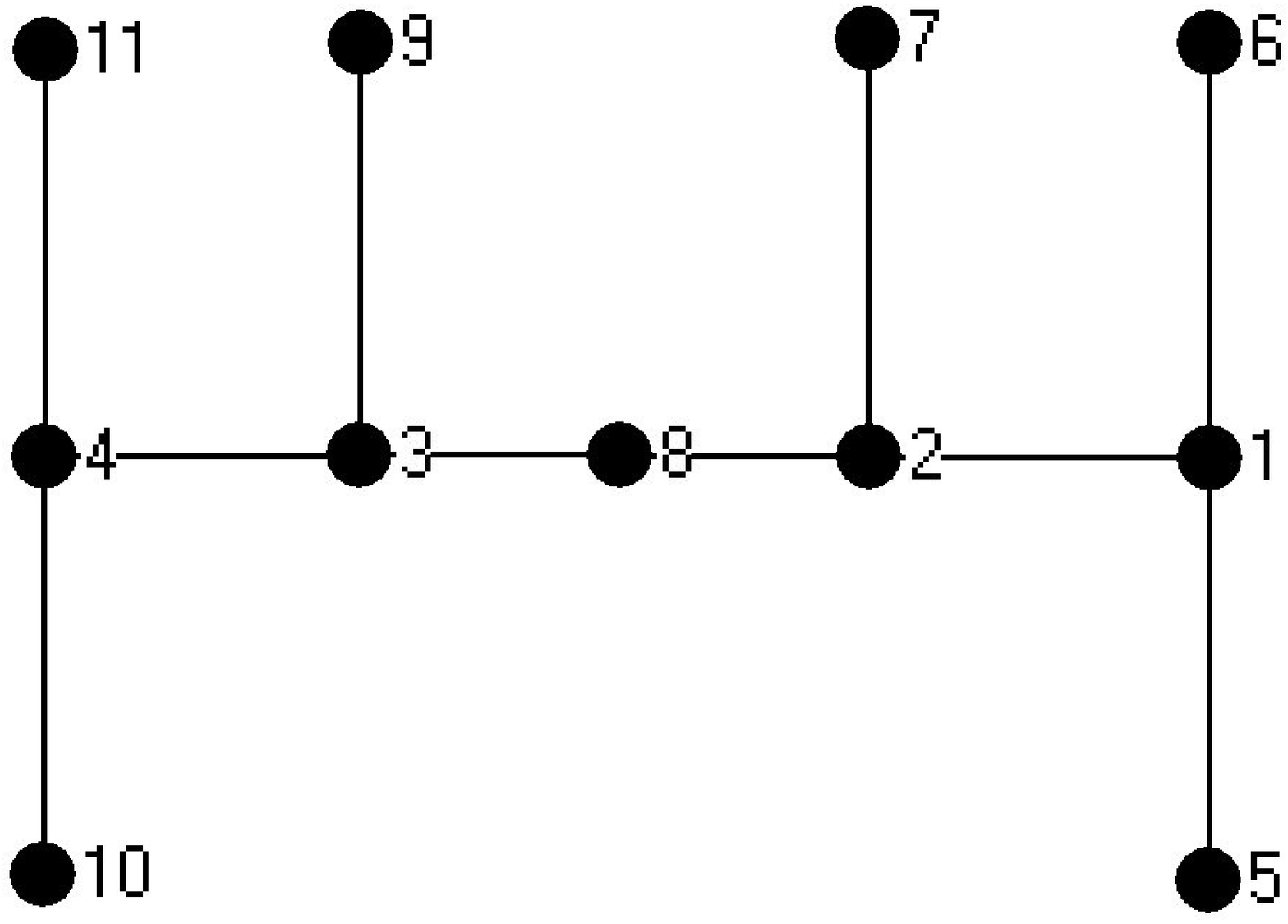}
\includegraphics[angle=90, width=6cm]{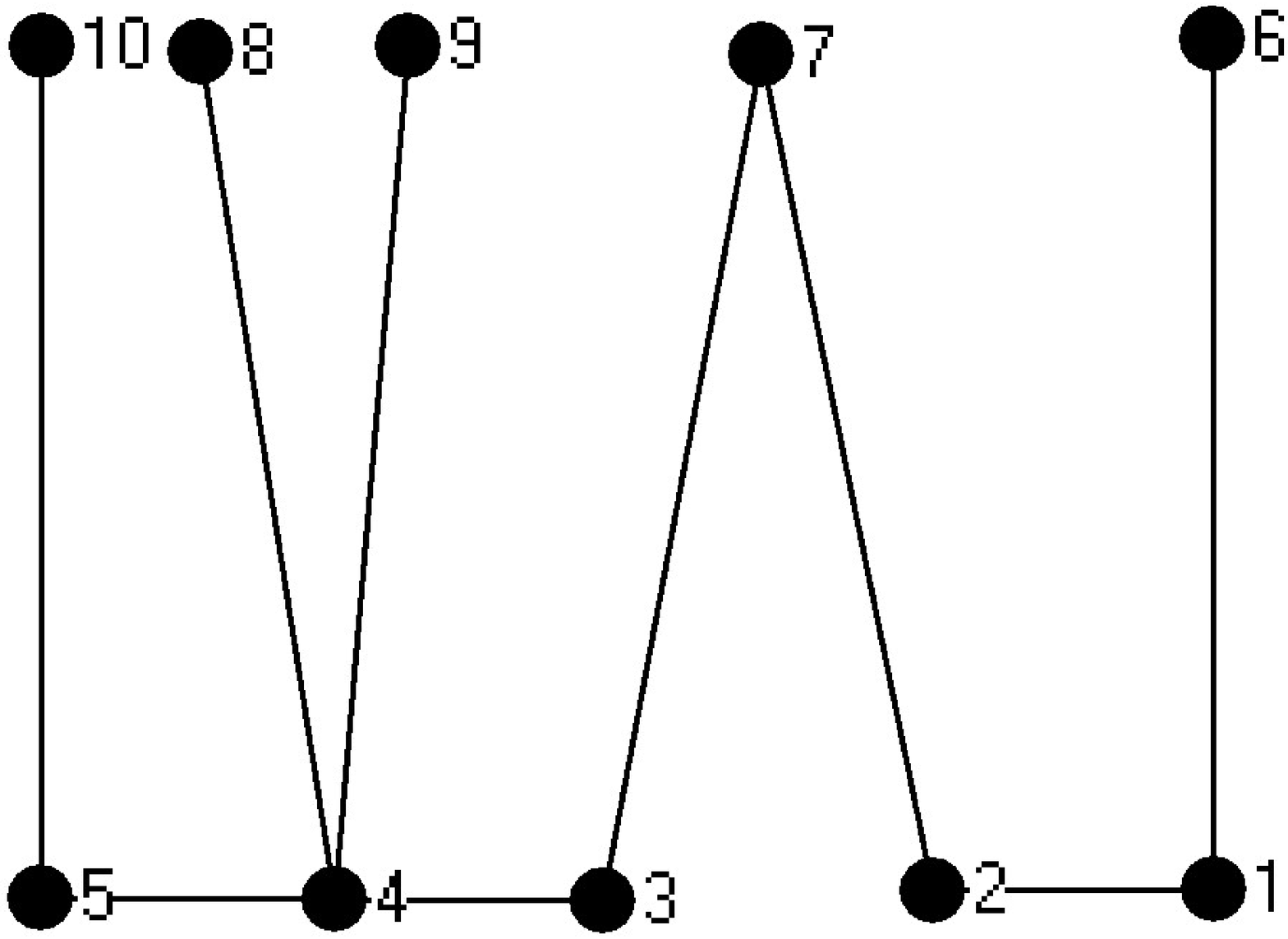}
\end{center}
\end{figure}

The above bound is attained, for instance,  for the left hand side
graph of Figure \ref{ej3}, where $S=\{1,2,3,4\}$ is a global
defensive $(-1)$-alliance and $\langle S\rangle$ has two connected
components. Moreover, the bound is attained in the case of the right
hand side graph of Figure \ref{ej3}, where $S=\{1,2,3,4,5\}$ is a
global defensive $0$-alliance and $\langle S\rangle$ has two
connected components.

\begin{corollary} For any tree $T$ of order $n$,
$\gamma_{k}^{a}(T)\ge
\left\lceil\displaystyle\frac{n+2}{3-k}\right\rceil.$
\end{corollary}

The above bound is attained for $k\in \{-4,-3,-2,0,1\}$ in the case
of $\Gamma=K_{1,4}$. As a particular case of above theorem  we
obtain the bounds obtained in \cite{GlobalalliancesOne}:
$$\gamma^{a}_{-1}(T)\ge \left\lceil\frac{n+2}{4}\right\rceil \quad
{\rm and }\quad \gamma^{a}_0(T)\ge
\left\lceil\frac{n+2}{3}\right\rceil.$$


\section{Global connected defensive  $k$-alliances}

It is clear that a defensive $k$-alliance of minimum cardinality
must induce a connected subgraph. But we can have a global defensive
$k$-alliance of minimum cardinality with nonconnected  induced
subgraph. We say that a defensive $k$-alliance $S$ is connected if
$\langle S\rangle$ is connected. We denote by $
\gamma_{k}^{ca}(\Gamma)$ the minimum cardinality of a global
connected defensive $k$-alliance in $\Gamma$. Obviously, $
\gamma_{k}^{ca}(\Gamma)\ge  \gamma_{k}^{a}(\Gamma)$. For instance,
for the left hand side graph of Figure \ref{ej3} we have
$\gamma_{-1}^{ca}(\Gamma)=5 > 4  = \gamma_{-1}^{a}(\Gamma)$ and for
the right hand side graph of Figure \ref{ej3} we have
$\gamma_{0}^{ca}(\Gamma)=6 > 5  = \gamma_{0}^{a}(\Gamma)$.

\begin{theorem} \label{diametADG}
For any connected graph $\Gamma$  of  diameter $D(\Gamma)$,
\begin{itemize}
\item[{\rm (}i{\rm )}]{$\gamma_{k}^{ca}(\Gamma)\ge
\left\lceil\frac{\sqrt{4(D(\Gamma)+n-1)+(1-k)^{2}}+(k-1)}{2}
\right\rceil.$} \item[{\rm (}ii{\rm )}]{$\gamma_{k}^{ca}(\Gamma)\ge
\left\lceil\frac{n+D(\Gamma)-1}{\left\lfloor\frac{\Delta-k}{2}\right\rfloor+2}\right\rceil
$.}
\end{itemize}
\end{theorem}
\begin{proof}
If $S$ is a dominating set in $\Gamma$ such that  $\langle S
\rangle$ is connected, then $D(\Gamma)\leq D(\langle S \rangle)+2.$
Hence,
\begin{equation} \label{cotadiam}
D(\Gamma)\leq |S|+1.
\end{equation}
Moreover, if $S$ is a global defensive  $k$-alliance  in $\Gamma$,
then $|S|$ satisfies (\ref{inecuacion}). The first result follows by
(\ref{inecuacion}) and (\ref{cotadiam}).

As a consequence of (\ref{cotavecinos3}), (\ref{strongGrado})  and
(\ref{cotadiam}) we obtain the second result.
\end{proof}

Both bounds  in Theorem \ref{diametADG} are tight. For instance,
both bounds are attained for $k\in \{-2,-1,0\}$ for the  graph of
Figure \ref{planarnue}. In such a case, both bounds lead to
$\gamma_{k}^{ca}(\Gamma)\geq 3.$  Moreover, both bounds lead to the
exact values of $\gamma_{k}^{ca}(K_{3,3})$  in the following cases:
$2\leq\gamma_{-3}^{ca}(K_{3,3})=\gamma_{-2}^{ca}(K_{3,3})=\gamma_{-1}^{ca}(K_{3,3}).$
Furthermore, notice  that bound (ii) leads to the exact values of
$\gamma_{k}^{ca}(Q_{3})$ in the cases
$4\leq\gamma_{0}^{ca}(Q_{3})=\gamma_{1}^{ca}(Q_{3})$, while bound
(i) only gives $3\leq\gamma_{0}^{a}(Q_{3})$ and $3\leq
\gamma_{1}^{ca}(Q_{3})$.

By Theorem \ref{diametADG}, and taking into a count that
$D(\Gamma)-1\le D({\cal L}(\Gamma))$, we obtain the following result
on the global connected $k$-alliance number of the line graph of
$\Gamma$ in terms of some parameters of $\Gamma$.

\begin{corollary}\label{agfConnected}
For any connected graph  $\Gamma$ of size $m$, diameter $D(\Gamma)$,
and maximum degrees $d_{1}\ge d_{2}$,
\begin{itemize}
\item[{\rm (}i{\rm )}] $\gamma_{k}^{ca}({\cal L}(\Gamma))\ge
\left\lceil\frac{\sqrt{4(D(\Gamma)+m-2)+(1-k)^{2}}-(1-k)}{2}
\right\rceil.$
\item[{\rm (}ii{\rm )}] $ \gamma_{k}^{ca}({\cal L}(\Gamma))\ge
\left\lceil\frac{2(m+D(\Gamma)-2)}{d_1+d_2-k+1}\right\rceil .$
\end{itemize}
\end{corollary}


\end{document}